\documentclass[reqno]{amsart}

\usepackage{amssymb,amscd,bbm,verbatim,graphicx,mathtools}
\usepackage[colorlinks]{hyperref}
\newtheorem{thm}{Theorem}[section]
\newtheorem{cor}[thm]{Corollary}
\newtheorem{lem}[thm]{Lemma}

\theoremstyle{definition}

\newtheorem{hyp}[thm]{Hypothesis}

\numberwithin{equation}{section}

\newcommand{\lVVert}{\left\vert\kern-0.25ex\left\vert\kern-0.25ex\left\vert}
\newcommand{\rVVert}{\right\vert\kern-0.25ex\right\vert\kern-0.25ex\right\vert}
\newcommand{\abs}[1]{\left\lvert#1\right\rvert}
\newcommand{\norm}[1]{\left\lVert #1\right\rVert}
\newcommand{\tnorm}[1]{\lVVert#1\rVVert}

\newcommand{\ov}{\overline}

\renewcommand{\Im}{\operatorname{Im}}
\newcommand{\supp}{\operatorname{supp}}

\newcommand{\dom}{\operatorname{dom}}
\newcommand{\ran}{\operatorname{ran}}
\newcommand{\rk}{\operatorname{rk}}

\newcommand{\BV}{\operatorname{BV}}
\newcommand{\BVl}{\operatorname{BV}_{\rm loc}}

\newcommand{\sgn}{\operatorname{sgn}}
\newcommand{\diag}{\operatorname{diag}}
\newcommand{\sm}[1]{\big(\begin{smallmatrix}#1\end{smallmatrix}\big)}
\newcommand{\bsm}[1]{\begin{pmatrix}#1\end{pmatrix}}
\newcommand{\bb}[1]{{\mathbb{#1}}}
\newcommand{\mc}[1]{{\mathcal{#1}}}

\newcommand{\id}{\mathbbm 1}

\newcommand{\ol}{{\overline{\lambda}}}
\newcommand{\la}{\lambda}

\newcommand{\mx}{{\rm max}}
\newcommand{\mn}{{\rm min}}
\newcommand{\fl}[2]{\lfloor#1\rfloor_{#2}}

\newcommand{\<}{\langle}
\renewcommand{\>}{\rangle}

\begin{document}
\title[Green's functions]{Green's functions for first-order systems of ordinary differential equations without the unique continuation property}
\author{Steven Redolfi and Rudi Weikard}
\address{Department of Mathematics, University of Alabama at Birmingham, Birmingham, AL 35226-1170, USA}
\email{stevenre@uab.edu, weikard@uab.edu}
\date{11. May 2022}
\thanks{This is a preprint of an article published in {\emph Integral Equations and Operator Theory} which is available online at \texttt{https://doi.org/10.1007/s00020-022-02703-6}.}
\thanks{\copyright 2022. This manuscript version is made available under the CC-BY-NC-ND 4.0 license \texttt{http://creativecommons.org/licenses/by-nc-nd/4.0/}}

\begin{abstract}
This paper is a contribution to the spectral theory associated with the differential equation $Ju'+qu=wf$ on the real interval $(a,b)$ when $J$ is a constant, invertible skew-Hermitian matrix and $q$ and $w$ are matrices whose entries are distributions of order zero with $q$ Hermitian and $w$ non-negative.
Under these hypotheses it may not be possible to uniquely continue a solution from one point to another, thus blunting the standard tools of spectral theory.
Despite this fact we are able to describe symmetric restrictions of the maximal relation associated with $Ju'+qu=wf$ and show the existence of Green's functions for self-adjoint relations even if unique continuation of solutions fails.
\end{abstract}
\maketitle

\section{Introduction}
This paper is a contribution to the spectral theory for the differential equation
$$Ju'+qu=wf$$
posed on the real interval $(a,b)$ when $J$ is a constant, invertible, and skew-Hermitian $n\times n$-matrix while the entries of the matrices $q$ and $w$ are distributions of order zero\footnote{Recall that distributions of order $0$ are distributional derivatives of functions of locally bounded variation and hence may be thought of, on compact subintervals of $(a,b)$, as measures. For simplicity we might use the word measure instead of distribution of order $0$ below.} with $q$ Hermitian and $w$ non-negative.
Ghatasheh and Weikard \cite{MR4047968} studied this equation under the additional hypothesis that initial value problems have unique balanced\footnote{A function of locally bounded variation is called balanced, if its values at any given point are averages of its left- and right-hand limits at that point.} solutions in the space of functions of locally bounded variation.

The equation $Ju'+qu=wf$ has, of course, been investigated by many people when the coefficients $q$ and $w$ are locally integrable.
In that situation initial value problems always have unique solutions.
This is not necessarily the case when the measures induced by $q$ or $w$ have discrete components.
It appears that an equation with measure coefficients was first considered in 1952, when Krein \cite{MR0054078} modelled a vibrating string.
In 1964 Atkinson \cite{MR0176141} suggested to unify the treatment of differential and difference equations by writing them as systems of integral equation where integrals were to be viewed as matrix-valued Riemann-Stieltjes integrals.
Atkinson explained that the presence of point masses may prevent the continuation of solutions across such points and posed a condition avoiding that problem but more restrictive than the one posed in \cite{MR4047968}.
In 1999 Savchuk and Shkalikov \cite{MR1756602} treated Schr\"odinger equations with potentials in the Sobolev space $W^{-1,2}_{\rm loc}$.
Their paper was very influential and spurred many further developments.
Nevertheless, Eckhardt et al. \cite{MR3046408} showed in 2013, with the help of quasi-derivatives or, equivalently, by writing the equation as a system, that a treatment without leaving the realm of locally integrable coefficients is possible.
In the same year Eckhardt and Teschl \cite{MR3095152} investigated $2\times 2$-systems with diagonal measure-valued matrices $q$ and $w$ requiring essentially Atkinson's condition.

A more thorough account of the subject's history is given in \cite{MR4047968}.
The papers \cite{MR3046408} and \cite{MR3095152}, mentioned above, may also serve as excellent sources, with perhaps different emphases, of this history.

One feature of systems of first-order equations is that, generally, they are represented by linear relations rather than linear operators.
There is a well-developed spectral theory for linear relations initiated by Arens \cite{MR0123188}, see also Orcutt \cite{Orcutt}, and Bennewitz \cite{MR0415380}.
The most important results (for our purposes) are also surveyed in Appendix B of \cite{MR4047968}.

Existence or uniqueness of solutions of an initial value problem for $Ju'+qu=wf$ fails when, for some $x\in(a,b)$, the matrices
$$B_\pm(x,0)=J\pm\frac12\Delta_q(x)$$
are not invertible.
Here $\Delta_q(x)=Q^+(x)-Q^-(x)$ when $Q$ denotes an anti-derivative of $q$.
Equivalently, $\Delta_q(x)=dQ(\{x\})$ where $dQ$ is the measure (locally) generated by $q$.
Assuming the unique continuation property for solutions of $Ju'+qu=wf$ Ghatasheh and Weikard defined maximal and minimal relations $T_\mx$ and $T_\mn$ associated with the differential equation $Ju'+qu=wf$ and showed that $T_\mx$ is the adjoint of $T_\mn$.
They characterized the self-adjoint restrictions of $T_\mx$, if any, with the aid of boundary conditions and proved that resolvents are given as integral operators, i.e., the existence of a Green's function for any such self-adjoint relation $T$.
Under even more restrictive conditions they also showed the existence of a Fourier transform diagonalizing $T$.

Campbell, Nguyen, and Weikard \cite{MR4298818} defined maximal and minimal relations and showed that $T_\mx=T_\mn^*$ without the hypothesis of unique continuation of solutions.
Our goal here is to advance their ideas.
In particular, even though the equation $Ju'+qu=w(\la u+f)$ may have infinitely many linearly independent solutions the deficiency indices, i.e., the number of linearly independent solutions of $Ju'+qu=\pm i wu$ of finite positive norm, is still bounded by $n$, the size of the system.
We show that symmetric restrictions of $T_\mx$, in particular the self-adjoint ones, are still given by posing boundary conditions and we show that the resolvents of self-adjoint restrictions are integral operators by proving the existence of Green's functions.

We will not approach the problem of Fourier transforms and eigenfunction expansions but hope to return to it in future work.

The material in this paper is arranged as follows.
In Section \ref{S.pre} we recall the circumstances under which existence and uniqueness of solutions to initial value problems does hold and investigate the sets of those $x\in(a,b)$ and $\la\in\bb C$ giving rise to trouble. Then, in Section \ref{S2} we discuss the manifold of solutions of our differential equation in the special case when $a$ and $b$ are regular endpoints.
These results are instrumental in Section \ref{S3} where we investigate the deficiency indices of the minimal relation and its symmetric extensions but without the assumption that $a$ and $b$ are regular.
Before we prove the existence of Green's functions for self-adjoint restrictions of the maximal relation in Section \ref{S:G} we discuss the role played by non-trivial solutions of zero norm in Section \ref{S:L0}.

Let us add a few words about notation.
$\mc D^{\prime0}((a,b))$ is the space of distributions of order $0$, i.e., the space of distributional derivatives of functions of locally bounded variation.
Any function $u$ of locally bounded variation has left- and right-hand limits denoted by $u^-$ and $u^+$, respectively.
Also, $u$ is called balanced if $u=u^\#=(u^++u^-)/2$.
The space of balanced functions of bounded variation defined on $(a,b)$ is denoted by $\BV^\#((a,b))$ while $\BVl^\#((a,b))$ stands for the space of balanced functions of \emph{locally} bounded variation.
We use $\id$ to denote an identity matrix of appropriate size and superscripts ${}^\top$ and ${}^*$ indicate transposition and adjoint, respectively.
The sum of two closed only trivially intersecting subspaces $S$ and $T$ of some Hilbert space (i.e., their direct sum) is denoted by $S\uplus T$; if $S$ and $T$ are even orthogonal we may use $\oplus$ instead of $\uplus$.
The orthogonal complement of a subspace $S$ of a Hilbert space $H$ is denoted by $H\ominus S$ or by $S^\perp$.
For $c_1, ..., c_N\in\bb C^n$ we abbreviate the column vector $(c_1^\top, ..., c_N^\top)^\top\in\bb C^{nN}$ by $(c_1, ..., c_N)^\diamond$.

\section{Preliminaries}\label{S.pre}
Throughout this paper we assume the following hypothesis to be in force.
\begin{hyp}\label{H:2.1}
$J$ is a constant, invertible and skew-Hermitian $n\times n$-matrix.
Both $q$ and $w$ are in $\mc D^{\prime0}((a,b))^{n\times n}$, $w$ is non-negative and $q$ Hermitian.
\end{hyp}

Given that $w$ is non-negative it gives rise to a positive measure on $(a,b)$ and we denote the space of functions $f$ which satisfy $\int f^*wf<\infty$ by $\mc L^2(w)$.
This space permits the semi-inner product $\<f,g\>=\int f^*wg$ (note that $\<f,f\>$ may be $0$ without $f$ being $0$).

Consider the differential equation
\begin{equation}\label{eq0}
Ju'+(q-\la w)u=wf
\end{equation}
where $\la$ is a complex parameter and $f$ an element of $\mc L^2(w)$.
The latter condition guarantees that $wf$ is in $\mc D^{\prime0}((a,b))^n$.
We will search for solutions in $\BVl^\#((a,b))^n$.
In this case each term in \eqref{eq0} is a distribution of order $0$ so that it makes sense to pose the equation.

The point $a$ is called a regular endpoint for $Ju'+qu=wf$, if there is a point $c\in(a,b)$ such that the left-continuous anti-derivatives $Q$ and $W$ of $q$ and $w$ are of bounded variation on $(a,c)$.
In this case $q$ and $w$ may be thought of as finite measures on $(a,c)$.
Similarly, $b$ is called regular, if $Q$ and $W$ are of bounded variation on $(c,b)$.
If an endpoint is not regular, it is called singular.
Not surprisingly, the study of our problem is less complicated when the endpoints are regular and we will use this fact to our advantage.

Despite our earlier denigration of the existence and uniqueness theorem of solutions of initial value problems it continues to play a crucial role.
The following theorem was proved in \cite{MR4047968}.
\begin{thm}\label{EUIVP}
Suppose $r\in \mc D^{\prime0}((a,b))^{n\times n}$, $g\in\mc D^{\prime0}((a,b))^{n}$ and that the matrices $\id\pm\Delta_r(x)/2$ are invertible for all $x\in(a,b)$.
Let $x_0$ be a point in $(a,b)$.
Then the initial value problem $u'=ru+g$, $u(x_0)=u_0\in\bb C^n$ has a unique balanced solution $u\in\BVl^\#((a,b))^n$.

If $a$ is a regular endpoint we may pose an initial condition (for $u^+$) at $a$.
Similarly, if $b$ is regular we may prescribe $u^-(b)$ as the initial condition.
\end{thm}

Suppose now that $u$ is a solution of \eqref{eq0}.
Treating either side of this equation as a measure (restricted to a compact subset of $(a,b)$) evaluation at a singleton $\{x\}$ shows that
$$J(u^+(x)-u^-(x))+\Delta_{q-\la w}(x)u^\#(x)=\Delta_w(x)f(x)$$
or, equivalently,
\begin{equation}\label{eq1}
B_+(x,\la)u^+(x)-B_-(x,\la)u^-(x)=\Delta_w(x)f(x)
\end{equation}
when we define
$$B_\pm(x,\la)=J\pm\frac12\big(\Delta_q(x)-\la\Delta_w(x)\big).$$
Note that, if $B_+(x,\la)$ is not invertible, we could be in one of the following two situations: (i) a solution given on $(a,x)$ may fail to exist on $(x,b)$ or (ii) there are infinitely many ways to continue a solution on $(a,x)$ to $(x,b)$.
An analogous statement holds, of course, if $B_-(x,\la)$ is not invertible.

Let us now investigate the circumstances when a pair $(x,\la)$ gives such trouble.
Define the sets $\Lambda_x=\{\la\in\bb C: \det(B_+(x,\la))\det(B_-(x,\la))=0\}$ and $\Xi_\la=\{x\in(a,b): \det(B_+(x,\la))\det(B_-(x,\la))=0\}$.
First note, since $B_-(x,\la)=-B_+(x,\ol)^*$, we have that $\Xi_\la=\Xi_\ol$ and that each $\Lambda_x$ is symmetric with respect to the real axis.
Also, $\Lambda_x$ is empty unless at least one of $\Delta_q(x)$ and $\Delta_w(x)$ is different from $0$ and hence for all but countably many $x$.
Next, we claim that $\Lambda_x$ is finite as soon as it misses one point.
To see this suppose that $B_+(x,\la_0)$ is invertible and that $\la\neq\la_0$.
Since
$$B_+(x,\la)=(\la_0-\la)B_+(x,\la_0)\big(\frac12 B_+(x,\la_0)^{-1}\Delta_w(x)-1/(\la-\la_0)\big)$$
we see that $B_+(x,\la)$ fails to be invertible only if $1/(\la-\la_0)$ is an eigenvalue of some $n\times n$-matrix.
A similar statement holds, of course, for $B_-$ proving our claim.

The really bad points $x$, namely those where $\Lambda_x=\bb C$, are thus contained in $\Xi_0$.
Here we wish to remove the hypothesis $\Xi_0=\emptyset$  posed in \cite{MR4047968}.
On any subinterval of $(a,b)$ on which $q$ gives rise to a finite measure we find that $\sum_{k=1}^\infty \|\Delta_q(x_k)\|$ must be finite,
when $k\mapsto x_k$ is a sequence of distinct points in that interval.
It follows now that $\Xi_0$ is a discrete set.
One shows similarly that, for any fixed complex number $\la$ the set $\Xi_\la$ is discrete.

\begin{lem}
Suppose $[s,t]\subset(a,b)$ and $(s,t)\cap \Xi_0=\emptyset$.
Then we have that $\Lambda_{(s,t)}=\bigcup_{x\in(s,t)} \Lambda_x$ is a discrete subset of $\bb C$.
\end{lem}

\begin{proof}
There are only finitely many points $x$ in $(s,t)$ where ${\|J^{-1}\Delta_q(x)\|>1}$.
Using a Neumann series one sees that only at such points the norm of $B_+(x,0)^{-1}$ can be larger than $2\|J^{-1}\|$.
Thus there is a positive number $C$ such that $\|B_+(x,0)^{-1}\|\leq C$ for all $x\in(s,t)$.
Now suppose that $B_+(x,\la)$ is not invertible and that $|\la|\leq R$.
Then $1/\la$ is an eigenvalue of $\frac12B_+(x,0)^{-1}\Delta_w(x)$.
This requires that $\|\Delta_w(x)\|\geq 2/(RC)$ and thus can happen only for finitely many $x\in(s,t)$.
Since similar arguments work for $B_-$ the number of points in $\bigcup_{x\in(s,t)} \Lambda_x$ which lie in a disk of radius $R$ centered at $0$ must be finite.
\end{proof}

We remark that, when one of the anti-derivatives of $q$ and $w$ is only locally of bounded variation, the set $\bigcup_{x\in(a,b)} \Lambda_x$ need not be discrete even if every $\Lambda_x$ is finite.

\begin{thm}\label{T:n2.3}
Suppose $[s,t]\subset(a,b)$ and $(s,t)\cap\Xi_0=\emptyset$.
If $u_0\in\bb C^n$ and $\la\in\bb C\setminus\Lambda_{(s,t)}$, then the initial value problem $Ju'+qu=\la wu$, $u^+(s)=u_0$ has a unique balanced solution in $(s,t)$.
Moreover, $u(x,\cdot)$ for $x\in(s,t)$ as well as $u^-(t,\cdot)$ are analytic in $\bb C\setminus\Lambda_{(s,t)}$ and meromorphic on $\bb C$.
An analogous statement holds when the initial condition is posed at $t$.
\end{thm}

\begin{proof}
The first claim is simply a consequence of Theorem \ref{EUIVP}.
When $x\in(s,t)$ the analyticity of $u(x,\cdot)$ in $\bb C\setminus\Lambda_{(s,t)}$, which is an open set, was proved in Section 2.3 of \cite{MR4047968}.
If we modify $q$ and $w$ by setting them $0$ on $[t,b)$ we do not change the solution on $(s,t)$.
The solution for the modified problem evaluated at $t$ is analytic and coincides with $u^-(t,\cdot)$ proving its analyticity.
It remains to show that a point $\la_0\in\Lambda_{(s,t)}$ can merely give rise to poles.

We know already that there are only finitely many points $x$ in $(s,t)$ where one of $B_\pm(x,\la_0)$ fails to be invertible.
Suppose $x'$ and $x''$ are two consecutive such points.
If we know the solution on $(s,x')$ and that $u^-(x',\cdot)$ has, at worst, a pole at $\la_0$, then the solution in $(x',x'')$ is determined by the initial value
$$u^+(x',\la)=B_+(x',\la)^{-1}B_-(x',\la)u^-(x',\la)$$
which also has, at worst, a pole at $\la_0$ since this is true for $B_+(x',\la)^{-1}$.
For $x\in(s,t)$ the claim follows now by induction.
To prove that $u^-(t,\cdot)$ is also meromorphic we proceed as before and modify $q$ and $w$ on $[t,b)$.
\end{proof}

\section{Solving the differential equation}\label{S2}
Our goal in this section is to investigate the set of solutions of the differential equation $Ju'+(q-\lambda w)u=wf$ on $(a,b)$ under a strengthened hypothesis.
\begin{hyp}\label{H:3.1}
In addition to Hypothesis \ref{H:2.1} we ask that $a$ and $b$ are regular endpoints for $Ju'+qu=wf$.

Moreover, given the partition
\begin{equation}\label{20200902.1}
a=x_0<x_1<x_2<...<x_N<x_{N+1}=b
\end{equation}
of $(a,b)$ we require that  $\Xi_0\subset\{x_1,...,x_N\}$.
We then consider only $\la$ for which both $B_+(x,\lambda)$ and $B_-(x,\la)$ are invertible unless $x$ is in $\{x_1,...,x_N\}$.
\end{hyp}
This hypothesis is in force throughout this section but later only if explicitly mentioned.
We emphasize that $\Xi_0$ is finite when $a$ and $b$ are regular.
Also, the set of permissible $\la$, which we call $\Omega_0$, is symmetric with respect to the real axis and avoids only a discrete set.

On each interval $(x_j,x_{j+1})$ we let $U_j(\cdot,\lambda)$ be a fundamental matrix of balanced solutions of the homogeneous differential equation $Ju'+(q-\lambda w)u=0$ such that $\lim_{x\downarrow x_j} U_j(x,\lambda)=\id$.
The existence of these fundamental matrices is guaranteed by Theorem \ref{EUIVP}.
The general balanced solution $u$ of the non-homogeneous equation $Ju'+(q-\la w)u=wf$ on $(x_j,x_{j+1})$ satisfies, according to Lemma 3.3 in \cite{MR4047968},
\[u^-(x)=U_j^-(x,\lambda)\big(c_j+J^{-1}\int_{(x_j,x)}U_j(\cdot,\ol)^*wf\big)\]
for any $c_j\in\bb C^n$.
Define \[U_j(x_{j+1},\lambda)=\lim_{x\uparrow x_{j+1}}U_j(x,\lambda)
\quad\text{and}\quad
I_j(f,\lambda)=\int_{(x_j,x_{j+1})}U_j(\cdot,\ol)^*wf.\]
Using $u^+(x_j)=c_j$ and $u^-(x_j)=U_{j-1}(x_j,\lambda)(c_{j-1}+J^{-1}I_{j-1}(f,\lambda))$ in equation \eqref{eq1} gives
\begin{multline*}
(-B_-(x_j,\lambda)U_{j-1}(x_j,\lambda),B_+(x_j,\la))\begin{pmatrix}c_{j-1}\\c_j\end{pmatrix}\\
 =\Delta_w(x_j)f(x_j)+B_-(x_j,\lambda)U_{j-1}(x_j,\lambda)J^{-1}I_{j-1}(f,\lambda).
\end{multline*}
We need to consider these equations for $j=1,...,N$ simultaneously.
This gives rise to the system
\begin{equation}\label{eq4}
\bb B(\lambda)\tilde u=\mc F_0(f,\la)
\end{equation}
where $\tilde u=(c_0,...,c_N)^\diamond$, $\bb B(\la)$, to be specified presently, is in $\bb C^{nN\times n(N+1)}$, and $\mc F_0(f,\la)$ is in $\bb C^{nN}$.
The two-diagonal block-matrix structure of $\bb B$ suggests the introduction of matrices $E_\top$ and $E_\bot$, which, respectively, strip the first and last $n$ components off a vector in their domain $\bb C^{n(N+1)}$.
If we also define the block-matrices
\[\mc B(\lambda)=\diag(B_+(x_1,\lambda),...,B_+(x_N,\lambda)),\]
\[\mc U(\lambda)=\diag(U_0(x_1,\lambda),...,U_{N-1}(x_N,\lambda)),\]
and $\mc J=\diag(J, ..., J)$
and when we note that
$$\mc B(\ol)^*=\diag(-B_-(x_1,\lambda),...,-B_-(x_N,\lambda)),$$
we obtain
\begin{equation}\label{bbB}
\bb B(\lambda)=\mc B(\ol)^*\mc U(\lambda)E_\bot+\mc B(\lambda)E_\top.
\end{equation}
The vector $\mc F_0(f,\la)$ is given by
$$\mc F_0(f,\la)=\mc R(f)-\mc B(\ol)^*\mc U(\lambda)\mc J^{-1}\mc I(f,\lambda)$$
with $\mc R(f)=((\Delta_wf)(x_1),...,(\Delta_wf)(x_N))^\diamond$ and $\mc I(f,\lambda)=(I_0(f,\lambda),...,I_{N-1}(f,\lambda))^\diamond$.

We now have the following theorem.
\begin{thm}\label{T:2.3}
The differential equation $Ju'+(q-\lambda w)u=wf$ has a solution $u$ on $(a,b)$ if and only if $\tilde u=(u^+(x_0),...,u^+(x_N))^\diamond$ is a solution of equation \eqref{eq4}.
In particular, in the homogeneous case, where $f=0$, the space of solutions has dimension ${n(N+1)}-\rk\bb B(\la)\geq n$.
\end{thm}
We note that $\rk\bb B(\la)=n$ when $N=1$ so that the space of solutions of $Ju'+(q-\la w)u=0$ is then exactly $n$-dimensional.
For $N=2$, however, consider the example $(a,b)=\bb R$, $J=\sm{0&-1\\ 1&0}$, $q=\sm{0&2\\ 2&0}(\delta_1-\delta_2)$, $w=\sm{2&0\\ 0&0}(\delta_1+\delta_2)$, where the $\delta_k$ are Dirac point measures concentrated on $\{k\}$.
It shows that the dimension of the space of solutions of $Ju'+(q-\la w)u=0$ may be strictly larger than $n$.

Next we investigate the connection between the right-hand limits of a solution $u$ of the homogeneous equation $Ju'+(q-\la w)u=0$ at the points $x_0$, ..., $x_N$ (given by the vector $\tilde u)$ and the vector $\hat u=(u(x_1), ..., u(x_N))^\diamond$.
We have $\hat u=\bb D(\la)\tilde u$ where
\begin{equation}\label{bbD}
\bb D(\la)=\frac{1}{2}(\mc U(\lambda)E_\bot+E_\top)
\end{equation}
is again a two-diagonal block-matrix.
If $N\geq 2$ we will also introduce the matrices $\bb B_m(\lambda)$ and $\bb D_m(\lambda)$ which are obtained by deleting the first and last $n$ columns from $\bb B(\lambda)$ and $\bb D(\lambda)$, respectively.
If $N=1$ we should think of $B_m(\lambda)$ and $D_m(\lambda)$ as maps from the trivial vector space to $\bb C^n$.
Their adjoints are the map from $\bb C^n$ to $\{0\}$.
With this understanding the following results hold also for $N=1$ even though they then involve ``matrices'' with no rows or columns.

\begin{lem} \label{L:2.4}
$\bb D(\ol)^*\bb B(\lambda)-\bb B(\ol)^*\bb D(\lambda)=\diag(-J,0,...,0,J)$ and $\bb D_m(\ol)^*\bb B(\lambda)-\bb B_m(\ol)^*\bb D(\lambda)=0$.
\end{lem}

\begin{proof}
This follows since $\mc U(\ol)^*\mc J \mc U(\la)=\mc J$ which, in turn, follows from Lemma 3.2 in \cite{MR4047968}.
\end{proof}

\begin{lem}\label{L:2.5}
The map $v\mapsto \bb B(\la)v$, restricted to $\ker\bb D(\la)$, is a bijection onto $\ker\bb D_m(\ol)^*$.
Similarly, the map $v\mapsto \bb D(\la)v$, restricted to $\ker\bb B(\la)$, is a bijection onto $\ker\bb B_m(\ol)^*$.
In particular, $\dim\ker\bb D(\la)=\dim\ker\bb D_m(\ol)^*$ and $\dim\ker\bb B(\la)=\dim\ker\bb B_m(\ol)^*$.
\end{lem}

\begin{proof}
The identity $\bb D_m(\ol)^*\bb B(\lambda)-\bb B_m(\ol)^*\bb D(\lambda)=0$ shows that $\bb B(\la)$ maps $\ker\bb D(\la)$ to $\ker\bb D_m(\ol)^*$ as well as that $\bb D(\la)$ maps $\ker\bb B(\la)$ to $\ker\bb B_m(\ol)^*$.

If $v\in\ker\bb B(\la)\cap\ker\bb D(\la)$ one shows that $E_\bot v=E_\top v=0$ using the definitions \eqref{bbB} and \eqref{bbD} of $\bb B$ and $\bb D$ and the fact that $\mc B(\la)-\mc B(\ol)^*=2\mc J$.
This, of course, implies that $v=0$ and hence the injectivity of both $\bb B(\la)|_{\ker\bb D(\la)}$ and $\bb D(\la)|_{\ker\bb B(\la)}$.

Clearly, both $\bb D(\la)$ and $\bb D_m(\ol)^*$, having invertible matrices along their main diagonal, are of full rank.
The rank-nullity theorem shows therefore that their kernels both have dimension $n$.
This proves surjectivity of $\bb B(\la)|_{\ker\bb D(\la)}$.

Finally, assume that $v\in\ker\bb B_m(\ol)^*$.
Then $v=\bb D(\la)x$ for some $x\in\bb C^{n(N+1)}$ which implies that $0=\bb B_m(\ol)^*\bb D(\la)x=\bb D_m(\ol)^*\bb B(\la)x$.
The first part of the proof shows that there is a $y\in\ker\bb D(\la)$ such that $\bb B(\la)y=\bb B(\la)x$.
Hence $v=\bb D(\la)(x-y)$ where $x-y\in\ker\bb B(\la)$.
\end{proof}

The following theorem establishes a connection between solutions of the differential equation $Ju'+(q-\la w)u=0$ and elements of $\ker\bb B_m(\ol)^*$.

\begin{thm} \label{T:2.6}
If $u$ is a solution of $Ju'+(q-\lambda w)u=0$ on $(a,b)$, then $\hat u=(u(x_1), ..., u(x_N))^\diamond$ is in $\ker\bb B_m(\ol)^*$.
If, in addition, $u^+(a)=u^-(b)=0$, then $\hat u\in\ker\bb B(\ol)^*$ (a subspace of $\ker\bb B_m(\ol)^*$).

Conversely, if $\hat u\in\ker\bb B_m(\ol)^*$, then $Ju'+(q-\lambda w)u=0$ has a unique solution $u$ on $(a,b)$ such that
$(u(x_1), ..., u(x_N))^\diamond=\hat u$.
If, indeed, $\hat u\in\ker\bb B(\ol)^*$, we further have $u^+(a)=u^-(b)=0$.
\end{thm}

Let us emphasize that $\supp u\subset[x_1,x_N]$ when $u^+(a)=u^-(b)=0$.

\begin{proof}
If $u$ solves $Ju'+(q-\lambda w)u=0$, then, by Theorem \ref{T:2.3}, $\tilde u\in\ker\bb B(\la)$.
Lemma~\ref{L:2.5} shows then that $\hat u=\bb D(\la)\tilde u$ is in $\ker\bb B_m(\ol)^*$.
If $u^+(a)=u^-(b)=0$, then Lemma~\ref{L:2.4} gives $0=\bb B(\ol)^*\bb D(\lambda)\tilde u=\bb B(\ol)^*\hat u$.

Conversely, assume that $\hat u\in\ker\bb B_m(\ol)^*=\bb D(\la)(\ker\bb B(\la))$.
Then there is a unique vector $\tilde u\in\ker\bb B(\la)$ such that $\hat u=\bb D(\la)\tilde u$, which, in turn, defines a unique solution $u$ of $Ju'+(q-\la w)u=0$ such that $(u(x_1), ..., u(x_N))^\diamond=\hat u$.
If $\hat u\in\ker\bb B(\ol)^*$, then, according to Lemma \ref{L:2.4}, $\diag(-J,0,...,0,J)\tilde u=0$ which shows that $u^+(a)=u^-(b)=0$.
\end{proof}

Given an algebraic system $Ax=b$ we know that there exist solutions only if $b\in\ran A=(\ker A^*)^\perp$.
For the differential equation $Ju'+(q-\la w)u=wf$ with integrable coefficients $q$ and $w$ the unique continuation property for the solutions gives rise to the variation of constants formula, which then guarantees the existence of solutions for any non-homogeneity $f$ (within reason).
In the present situation, however, the problem of existence raises its head and we now set out to give necessary and sufficient conditions for $f$ guaranteeing the existence of a solution in the spirit of Linear Algebra.

\begin{lem}\label{L:2.7}
If $\tilde v\in\ker\bb B(\ol)$ and $\hat v=\bb D(\ol)\tilde v$, then
$$\tilde v^*E_\bot^*=-\hat v^*\mc B(\ol)^*\mc U(\la)\mc J^{-1}\quad\text{and}\quad \tilde v^*E_\top^*=\hat v^*\mc B(\la)\mc J^{-1}.$$

Moreover, if $f\in\mc L^2(w)$ and $Jv'+(q-\ol w)v=0$, then
$$\int v^*wf=\hat v^*\mc F_0(f,\la)+\hat v^*\mc B(\la)\mc J^{-1}\tilde{\mc I}(f,\la)
=\hat v^*\mc F_0(f,\la)+\tilde v^*E_\top^*\tilde{\mc I}(f,\la)$$
where $(v(x_1),...,v(x_N))^\diamond=\hat v=\bb D(\ol)\tilde v$ and $\tilde{\mc I}(f,\la)=(0,...,0,I_N(f,\la))^\diamond\in\bb C^{nN}$.
\end{lem}

\begin{proof}
Using the definitions \eqref{bbB} and \eqref{bbD} of $\bb B$ and $\bb D$ and the identities $\mc B(\la)-\mc B(\ol)^*=2\mc J$ and
$\mc U(\la)^*\mc J \mc U(\ol)=\mc J$ we obtain that $\bb B(\ol)\tilde v=0$ implies
$$\mc B(\ol)\bb D(\ol)\tilde v=\mc U(\la)^{*-1}\mc JE_\bot\tilde v\quad\text{and}\quad
 \mc B(\lambda)^*\bb D(\ol)\tilde v=-\mc J E_\top\tilde v.$$
Taking adjoints gives the first claim since $\hat v=\bb D(\ol)\tilde v$.

The second claim is an immediate consequence of this, since
\begin{align*}
\int v^*wf
 &= \hat v^*\mc R(f)+ \tilde v^*(I_0(f,\la), ..., I_N(f,\la))^\diamond\\
 &= \hat v^*\mc R(f)+ \tilde v^*E_\bot^* \mc I(f,\la)+\tilde v^*E_\top^*\tilde{\mc I}(f,\la)\\
 &= \hat v^*\mc R(f)- \hat v^*\mc B(\ol)^*\mc U(\la)\mc J^{-1}\mc I(f,\la)+\hat v^*\mc B(\la)\mc J^{-1}\tilde{\mc I}(f,\la)\\
 &= \hat v^*\mc F_0(f,\la)+\hat v^*\mc B(\la)\mc J^{-1}\tilde{\mc I}(f,\la).
\end{align*}
\end{proof}

\begin{thm}\label{T:2.8}
The differential equation $Ju'+(q-\la w)u=wf$ has a solution on $(a,b)$ if and only if $\int v^*wf=0$ for every solution $v$ of $Jv'+(q-\ol w)v=0$ which vanishes at $a$ and $b$.
\end{thm}

\begin{proof}
By Theorem \ref{T:2.3} the solution $u$ exists if and only if the system \eqref{eq4} has a solution $\tilde u=(u^+(x_0),...,u^+(x_N))^\diamond$.
This, in turn, happens if and only if $\mc F_0(f,\la)\in\ran\bb B(\la)=(\ker\bb B(\la)^*)^\perp$.

By Theorem \ref{T:2.6} the solutions of $Jv'+(q-\ol w)v=0$ which vanish at $a$ and $b$ are in one-to-one correspondence with elements of $\ker\bb B(\la)^*$.
Since $v^+(x_N)=0$ we have $\tilde v^*E_\top^*\tilde{\mc I}(f,\la)=0$ and then, from Lemma \ref{L:2.7}, we obtain $\hat v^*\mc F_0(f,\la)=\int v^*wf$.
\end{proof}

In the case of unique continuation of solutions the condition that $v$ vanishes at $a$ or $b$ implies, of course, that $v=0$.
Consequently, $Ju'+{(q-\la w)u}=wf$ has then a solution for any $f\in\mc L^2(w)$.
The set of all solutions is thus obtained by adding the general solution of $Ju'+(q-\la w)u=0$ whose dimension is ${n(N+1)}-\rk\bb B(\la)\geq n$.

\begin{thm}
The differential equation $Ju'+(q-\la w)u=wf$ has a solution on $(a,b)$ which vanishes at $a$ and $b$ if and only if $\int v^*wf=0$ for every solution $v$ of $Jv'+(q-\ol w)v=0$.
\end{thm}

\begin{proof}
For $u$ to vanish at $a$ and $b$ it is required that $u^+(x_0)=0$ and $u^+(x_N)=-J^{-1}I_N(f,\la)$.
The system \eqref{eq4} is therefore equivalent to
$$\bb B_m(\la)(c_1,...,c_{N-1})^\diamond = \mc F_0(f,\la)+\mc B(\la)\mc J^{-1}\tilde{\mc I}(f,\la).$$
The proof is now analogous to the one for Theorem \ref{T:2.8}.
\end{proof}

We conclude this section by ``counting'' the solutions of $Ju'+qu=\la wu$ which are not compactly supported.
More precisely, we will determine the dimension of the quotient space of all solutions of $Ju'+qu=\la wu$ modulo the space of compactly supported solutions.
Theorem \ref{T:2.6} shows that the space of all solutions of $Ju'+qu=\la wu$ is in one-to-one correspondence with $\ker\bb B_m(\ol)^*$ and that the space of compactly supported solutions of $Ju'+qu=\la wu$ is in one-to-one correspondence with $\ker\bb B(\ol)^*$.
We therefore define
$$\tilde n(\la)=\dim(\ker\bb B_m(\ol)^*/\ker\bb B(\ol)^*)=\dim\ker\bb B_m(\ol)^*- \dim\ker\bb B(\ol)^*.$$

\begin{lem}
$\tilde n(\la)+\tilde n(\ol)=2n$.
\end{lem}

\begin{proof}
Since $\rk\bb B(\la)=\rk\bb B(\la)^*$, the rank-nullity theorem implies
$$\dim\ker\bb B(\la)=n(N+1) -\rk\bb B(\la)^*=n+\dim\ker\bb B(\la)^*.$$
Hence, using also the analogous equation for $\ol$,
$$\dim\ker\bb B(\la)-\dim\ker\bb B(\la)^*+\dim\ker\bb B(\ol)-\dim\ker\bb B(\ol)^*=2n.$$
Lemma \ref{L:2.5} gives that $\dim\ker\bb B(\la)=\dim\ker\bb B_m(\ol)^*$ yielding the claim.
\end{proof}

From Theorem \ref{T:n2.3} we know that the matrices $U_j(x_{j+1},\cdot)$ are meromorphic on $\bb C$ with poles at most at points in the complement of $\Omega_0$.
It follows that the entries of $\bb B$ are also meromorphic.
Since the meromorphic functions on $\bb C$ form a field there is a row-echelon matrix $\tilde{\bb B}$ with meromorphic entries such that $\bb B\tilde u=0$ has the same solutions as $\tilde{\bb B}\tilde u=0$.
Now define a set $\Omega$ as $\Omega_0$ without the set of all poles of $\tilde{\bb B}$ as well as their complex conjugates, and the set of zeros and their conjugates of any of the pivots of $\tilde{\bb B}$.

\begin{thm}\label{L:3.2}
If $\la\in\Omega$, then $\dim\ker\bb B(\la)=\dim\ker\bb B(\ol)$ and $\tilde n(\la)=n$.
\end{thm}

\begin{proof}
The construction of $\Omega$ entails that $\rk\bb B(\la)=\rk\tilde{\bb B}(\la)=\rk\bb B(\ol)$ if $\la\in\Omega$.
Since $\tilde n(\la)=\dim\ker\bb B(\la)-\dim\ker\bb B(\ol)^*=\dim\ker\bb B(\la)+n-\dim\ker\bb B(\ol)$ we obtain $\tilde n(\la)=n$.
\end{proof}

\section{Symmetric restrictions of \texorpdfstring{$T_\mx$}{Tmax}}\label{S3}
Given a differential equation $Ju'+qu=wf$ we now define associated minimal and maximal relations.
Recall that $\mc L^2(w)$ is the space of functions $f$ such that $\int f^*wf<\infty$.
First we define
$$\mc T_\mx=\{(u,f)\in \mc L^2(w)\times\mc L^2(w): u\in\BVl^\#((a,b))^n, Ju'+qu=wf\}.$$
Subsequently we will always tacitly assume that $u\in\BVl^\#((a,b))^n$, when we use $u'$.
Next, let
$$\mc T_\mn=\{(u,f)\in \mc T_\mx: \text{$\supp u$ is compact in $(a,b)$}\}.$$
Note that these are spaces of pairs of functions.
To employ the power of functional analysis we need to realize these relations in Hilbert spaces.
Therefore we introduce, as usual, the space $L^2(w)$ as the quotient of $\mc L^2(w)$ modulo the subspace of all $u\in\mc L^2(w)$ for which $\norm{u}^2=\int u^*wu=0$.
Denoting the equivalence class corresponding to $u$ by $[u]$ we now set
$$T_\mx=\{([u],[f])\in L^2(w)\times L^2(w): (u,f)\in \mc T_\mx\}$$
and
$$T_\mn=\{([u],[f])\in T_\mx: (u,f)\in \mc T_\mn\}.$$
Here (and elsewhere) we choose brevity over precision: whenever we have a pair $([u],[f])$ in $T_\mx$ we choose $u$ and $f$ such that $(u,f)\in \mc T_\mx$.

Define the vector space
\[\mathcal{L}_0=\{u\in\BVl^\#((a,b))^n: Ju'+qu=0\text{ and } \|u\|=0\}.\]
In many cases this space is trivial and some authors restrict their attention to the case where it is; this is then called the definiteness condition.
However, we will not do so here.
Note that $\|u\|=0$ if and only if $wu$ is the zero distribution.
The significance of $\mc L_0$ stems from the following fact.
Suppose $([u],[f])\in T_\mx$ and that there are $u,v\in[u]$ and $f,g\in[f]$ such that $Ju'+qu=wf$ and $Jv'+qv=wg$.
Then $J(u-v)'+q(u-v)=w(f-g)=0$ as well as $w(u-v)=0$, i.e., $u-v\in\mc L_0$.
In other words, in the presence of a non-trivial space $\mc L_0$, the class $[u]$ has many representatives of locally bounded variation satisfying the differential equation for a given class $[f]$ (the choice of a representative of $[f]$, on the other hand, is irrelevant).
In Section \ref{S:L0} we will describe a procedure to choose a representative of $[u]$ in a distinctive way.

In \cite{MR4298818} it was proved that $T_\mn$ is symmetric, indeed that $T_\mn^*=T_\mx$.
In this case it is well-known that von Neumann's theorem holds.
Setting $D_\la=\{([u],\la[u])\in T_\mx\}$ it states that
$$T_\mx=\ov{T_\mn}\uplus D_\la\uplus D_{\ol}$$
when $\Im\la\neq0$.
Moreover, when $\la=\pm i$, these direct sums are even orthogonal.
It is also known that the dimension of $D_\la$ does not change as $\la$ varies in either the upper or the lower half plane.
The numbers $n_\pm=\dim D_{\pm i}$ are called deficiency indices of $T_\mn$ and we are now setting out to investigate these.

If $u$ is a solution of $Ju'+qu=\la wu$ which is compactly supported then $(u,\la u)\in \mc T_\mn$ and $([u],\la[u])\in T_\mn\cap D_\la$.
If $\la$ is not real, then $T_\mn\cap D_\la$ is trivial and it follows that compactly supported solutions of $Ju'+qu=\la wu$ do not contribute to the corresponding deficiency index.
We now have, as a corollary of Theorem \ref{L:3.2}, that the deficiency indices of $T_\mn$ cannot be more than $n$ if $a$ and $b$ are regular endpoints.
We do not state this result separately since it is included in the next theorem about the general case.

Thus, to emphasize, we allow in the following $a$ and $b$ to be either regular or singular endpoints.
Let $\tau_k$, $k\in\bb Z$, be a strictly increasing sequence in $(a,b)$ having $a$ and $b$ as its only limit points and such that all points in $\Xi_0$ are among the $\tau_k$.
Considering now only the interval $I_k=(\tau_{-k},\tau_k)$ we set $x_j=\tau_{-k+j}$ for $j=0, ..., N+1=2k$.
We can then introduce the objects from Section \ref{S2}.
To emphasize their dependence on $k$ we will add a superscript $(k)$ to those objects.
We have then, in particular, the matrices $\bb B^{(k)}$, $\bb B_m^{(k)}$ and the sets $\Omega^{(k)}$ of permissible values of $\la$.
We now define $\Omega=\bigcap_{k=1}^\infty \Omega^{(k)}$ and note that $\Omega$ is symmetric with respect to the real axis and misses only countably many values from $\bb C$.

Now fix a non-real $\la\in\Omega$.
If $u$ is a solution of $Ju'+qu=\la wu$ on $(a,b)$ we denote its restriction to the interval $I_k$ by $u^{(k)}$.
We are interested in the quotient space $X_k$ of all solutions of $Ju'+qu=\la wu$ on $I_k$ modulo the compactly supported solutions.
If $u$ is a solution of $Ju'+qu=\la wu$ on $I_k$ we denote the associated equivalence class in $X_k$ by $\fl{u}{k}$.
A compactly supported solution $u$ of $Ju'+qu=\la wu$ on $I_k$ can be extended by $0$ to all of $(a,b)$ yielding an element in $T_\mn\cap D_\la$.
This implies, since $\Im\la\neq0$, that $\|u\|^2=\int_{I_k} u^*wu=0$ and shows that $X_k$ is a normed space with the norm given by $\|u\|_k^2=\int_{I_k} u^*wu$.
According to Theorem \ref{T:2.6} the quotient space $X_k$ is isomorphic to $\ker\bb B_m^{(k)}(\ol)^*/\ker\bb B^{(k)}(\ol)^*$ and, by Theorem \ref{L:3.2}, its dimension is equal to $n$ since $\la\in\Omega\subset\Omega^{(k)}$.

\begin{thm}
The deficiency indices of $T_\mn$ are less than or equal to $n$.
\end{thm}

\begin{proof}
Fix a non-real $\la\in\Omega$.
Suppose $u_1$, ..., $u_m$ are solutions of $Ju'+qu=\la wu$ such that $[u_1]$, ..., $[u_m]$ are linearly independent elements of $D_\la$.
We will show below that there is an interval $I_p=(\tau_{-p},\tau_p)$ such that $\fl{u_1^{(p)}}{p}$, ..., $\fl{u_m^{(p)}}{p}$ are linearly independent elements of $X_p$.
Hence $m\leq n$, the dimension of $X_p$.
Since deficiency indices are constant in either half-plane they cannot be larger than $n$.

We will now prove the existence of $I_p$ by induction.
That is we prove that, for every $k\in\{1,..., m\}$, there is an interval $I_{\ell_k}$ such that the restrictions of $u_1$, ..., $u_k$ to $I_{\ell_k}$ generate linearly independent elements $\fl{u_1^{(\ell_k)}}{\ell_k}$, ..., $\fl{u_k^{(\ell_k)}}{\ell_k}$ of $X_{\ell_k}$.
Once this is achieved we set $p=\ell_m$.

Suppose $k=1$ and let $I_{\ell_1}$ be an interval such that $\|u_1^{(\ell_1)}\|>0$.
By what we argued above we know that $u_1^{(\ell_1)}$ is not compactly supported in $I_{\ell_1}$ and thus gives rise to a non-zero (and hence linearly independent) element of $X_{\ell_1}$.

Now suppose we had already shown our claim for some $k<m$.
If $\fl{u_1^{(\ell_k)}}{\ell_k}$, ..., $\fl{u_{k+1}^{(\ell_k)}}{\ell_k}$ are already linearly independent as elements of $X_{\ell_k}$ we choose $\ell_{k+1}=\ell_k$ and our induction step is complete.
Otherwise, there are unique complex numbers $\alpha_1$, ..., $\alpha_k$ such that
$$\|(\alpha_1u_1+ ... +\alpha_k u_k+u_{k+1})^{(\ell_k)}\|_{\ell_k}=0.$$
However, there must be an interval $I_{\ell_{k+1}}\supset I_{\ell_k}$ where
$$\|(\alpha_1u_1+ ... +\alpha_k u_k+u_{k+1})^{(\ell_{k+1})}\|_{\ell_{k+1}}>0$$
on account that $[u_1]$, ..., $[u_{k+1}]$ are linearly independent.
It follows now that, as elements of $X_{\ell_{k+1}}$ the vectors $\fl{u_1^{(\ell_{k+1})}}{\ell_{k+1}}$, ..., $\fl{u_{k+1}^{(\ell_{k+1})}}{\ell_{k+1}}$ are linearly independent.
This completes our induction step also in this case.
\end{proof}

\begin{cor}
If $a$ and $b$ are regular, then $n_+=n_-$.
\end{cor}

\begin{proof}
Fix a non-real $\la$ in $\Omega$.
Since $a$ and $b$ are regular, the set $\Xi_\la=\Xi_\ol$ is finite.
Thus we may assume that it is contained in $I_k=(\tau_{-k},\tau_k)$ for some appropriate $k$.
Then $\dim \ker\bb B^{(k)}(\la)$ is the number of linearly independent solutions of $Ju'+qu=\la wu$.
Theorem \ref{L:3.2} shows that $Ju'+qu=\ol wu$ has the same number of linearly independent solutions.
Any of these solutions has finite norm but some may have norm $0$.
Now note, that if $u$ is a solution of $Ju'+qu=\la wu$ of norm $0$, then we have $wu=0$, so that $u$ is also a solution of $Ju'+qu=\ol wu$.
Therefore $n_+=n_-$.
\end{proof}

As mentioned above, it is well-known, even in the case of relations, that von Neumann's theorem $E^*=E\oplus D_i\oplus D_{-i}$ holds when $E$ is a closed symmetric relation in $\mc H\times\mc H$ when $\mc H$ is a Hilbert space.
In our case, when $d=\dim D_i\oplus D_{-i}$ is finite, as we just showed, we can use Theorem~B.5 in \cite{MR4047968} to characterize the symmetric restriction of $T_\mx$ in terms of boundary conditions.
We state that theorem here for easy reference.
The operator $\mc J$ appearing there is defined by $\mc J(u,f)=(f,-u)$ for $u,f\in\mc H$.

\begin{thm}\label{T:B.6}
Suppose $E$ is a closed symmetric relation in $\mc H\times\mc H$ with $d=\dim D_i\oplus D_{-i}<\infty$ and that $m\leq d/2$ is a natural number or $0$.
If $A:E^*\to\bb C^{d-m}$ is a surjective linear operator such that $E\subset \ker A$ and $A\mc JA^*$ has rank $d-2m$ then $\ker A$ is a closed symmetric restriction of $E^*$ for which the dimension of $(\ker A)\ominus E$ is $m$.
Conversely, every closed symmetric restriction of $E^*$ is the kernel of such a linear operator $A$.
Finally, $\ker A$ is self-adjoint if and only if $A\mc JA^*=0$ (entailing $m=d/2$).
\end{thm}

A second ingredient for our next considerations is Lagrange's identity (or Green's formula).
If $(u,f)$ and $(v,g)$ are in $\mc T_\mx$, then $v^*wf$ and $g^*wu$ are finite measures.
Therefore $v^*Ju'+v^{\prime*}Ju=v^*wf-g^*wu$ is also a finite measure.
Its antiderivative $v^*Ju$ is of bounded variation and thus has limits at $a$ and $b$.
Integration now gives Lagrange's identity
\begin{equation}\label{Lagrange}
(v^*Ju)^-(b)-(v^*Ju)^+(a)=\<v,f\>-\<g,u\>.
\end{equation}
Note the right-hand side, and hence the left-hand side, does not change upon choosing different representatives in place of $u,f,v$, or $g$.

Now, if $(v,g)$ is an element of $D_i\oplus D_{-i}$, then $(u,f)\mapsto \<(v,g),(u,f)\>$ is a bounded linear functional on $T_\mx$.
Conversely, since $T_\mx$ is a Hilbert space, a bounded linear functional on $T_\mx$ is given by $(u,f)\mapsto \<(v,g),(u,f)\>$ for some $(v,g)\in T_\mx$.
When it is also known that $\ov{T_\mn}$ is in the kernel of this functional, $(v,g)$ may be chosen in $D_i\oplus D_{-i}$.
Hence, in our situation, the operator $A$ from Theorem \ref{T:B.6} is given by $d-m$ linearly independent elements in $D_i\oplus D_{-i}$.
Lagrange's identity implies that the entries of the matrix $A\mc JA^*$ are then given by
\begin{equation}\label{200820.1}
(A\mc JA^*)_{k,\ell}=\<(v_k,g_k),(g_\ell,-v_\ell)\>=(g_k^*Jg_\ell)^-(b)-(g_k^*Jg_\ell)^+(a).
\end{equation}

Therefore we arrive at the following theorem.
\begin{thm}\label{T:3.5}
Let $d=n_++n_-$ and suppose that $m\leq\min\{n_+,n_-\}$.
If $(v_1,g_1)$, ..., $(v_{d-m},g_{d-m})$ are linearly independent elements of $D_i\oplus D_{-i}$ such that the matrix defined in \eqref{200820.1} has rank $d-2m$, then
\begin{equation}\label{200821.1}
T=\{(u,f)\in T_\mx: (g_j^*Ju)^-(b)-(g_j^*Ju)^+(a)=0\text{ for $j=1,...,d-m$}\}
\end{equation}
is a closed symmetric restriction of $T_\mx$.

Conversely, if $T$ is a closed symmetric restriction of $T_\mx$ and $m$ is the dimension of $T\ominus \ov{T_\mn}$, then $T$ is given by \eqref{200821.1} for appropriate elements $(v_1,g_1)$, ..., $(v_{d-m},g_{d-m})$ of $D_i\oplus D_{-i}$ for which the matrix defined in \eqref{200820.1} has rank $d-2m$.
\end{thm}

For self-adjoint restrictions of $T_\mx$ it is hence necessary and sufficient that $n_+=n_-=m=d-m$ and that $(g_k^*Jg_\ell)^-(b)-(g_k^*Jg_\ell)^+(a)=0$ for all $1\leq k,\ell\leq m=d/2$.

\section{The space \texorpdfstring{$\mc L_0$}{L0}}\label{S:L0}
We mentioned earlier that the class $[u]$ does not have a unique balanced representative when $([u],[f])\in T_\mx$, if the space $\mc L_0$ has non-trivial elements.
In this section we describe a procedure to choose a representative in a distinctive way.

To this end we assume, without loss of generality, that $B_+(\tau_0,0)=B_-(\tau_0,0)=J$ so that solutions of our differential equations are continuous at $\tau_0$.
Define $N_0=\{h(\tau_0):h\in\mc L_0\}$ and for each $k\in\bb N$ both $N_k=\{h^+(\tau_k): h\in\mc L_0, \supp h\subset[\tau_k,b)\}$ and
$N_{-k}=\{h^-(\tau_{-k}): h\in\mc L_0, \supp h\subset(a,\tau_{-k}]\}$.
Then, for $k\in\bb N_0$, we say that a function $u\in\BVl^\#((a,b))^n$ satisfies condition $(\pm k)$, if $u^\pm(\tau_{\pm k})$ is perpendicular to $N_{\pm k}$ (using always the upper sign or always the lower sign).

\begin{lem}\label{L:5.1}
Suppose $([u],[f])\in T_\mx$.
Then there is a unique balanced $v\in[u]$ such that $(v,f)\in\mc T_\mx$ and $v$ satisfies condition $(k)$ for every $k\in\bb Z$.
\end{lem}

\begin{proof}
First consider uniqueness.
Suppose $u$ and $v$ are two functions satisfying the given conditions.
Then $u-v\in\mc L_0$ and hence $(u-v)(\tau_0)^*t(\tau_0)=0 $ for $t=u$ and $t=v$.
Subtract these equations to find $(u-v)(\tau_0)=0$, and thus $u=v$ on $(\tau_{-1},\tau_1)$.
Moreover, $h_1=(u-v)\chi_{[\tau_1,b)}$ and $h_{-1}=(u-v)\chi_{(a,\tau_{-1}]}$ are in $\mc L_0$.
Conditions (1) and $(-1)$ show therefore that $(u-v)^+(\tau_1)$ and $(u-v)^-(\tau_{-1})$ are also $0$ which proves that $u=v$ on $(\tau_{-2},\tau_2)$.
Induction informs us now that $u=v$ everywhere.

We now turn to existence.
Pick a balanced representative $u\in[u]$ such that $(u,f)\in\mc T_\mx$.
There is an element $h_0\in\mc L_0$ such that the orthogonal projection of $u(\tau_0)$ onto $N_0$ equals $h_0(\tau_0)$.
Thus $v_0=u-h_0$ satisfies $(v_0,f)\in\mc T_\mx$, $v_0\in[u]$, and condition $(0)$.

Next, there is an element $h_1\in\mc L_0$ with support in $[\tau_1,b)$ such that the orthogonal projection of $v_0^+(\tau_1)$ onto $N_1$ equals $h_1^+(\tau_1)$.
We now define $v_1=v_0-h_1$.
Then $(v_1,f)\in\mc T_\mx$, $v_1\in[u]$, and $v_1$ satisfies condition $(1)$.
Notice that $v_1=v_0$ on $(a,\tau_1)$ implying that $v_1$ also satisfies condition $(0)$.

Proceeding recursively, we may define, for each $k\in\bb N$, functions $h_k\in\mc L_0$ supported in $[\tau_k,b)$ such that
$v_k=u-\sum_{j=0}^k h_j$ satisfies conditions $(0)$, ..., $(k)$, $v_k\in[u]$, and $(v_k,f)\in\mc T_\mx$.

Since, for a fixed $x\in[\tau_0,b)$, only finitely many of the numbers $h_k(x)$ are different from $0$, we find that the sequence $k\mapsto v_k$ converges pointwise to a function $\tilde v\in[u]$ satisfying conditions $(k)$ for all $k\in\bb N_0$ and $(\tilde v,f)\in\mc T_\mx$.
We can now repeat this process for negative integers starting from the function $\tilde v$ instead of $u$ arriving eventually at a function $v\in[u]$ satisfying conditions $(k)$ for all $k\in\bb Z$ and $(v,f)\in\mc T_\mx$.
\end{proof}

We denote the operator which assigns the function $v$ just constructed to a given element $([u],[f])\in T_\mx$ by $E$.
If $I_m=(\tau_{-m},\tau_m)$ we also define $E_m:T_\mx\to\BV^\#(I_m)^n$ by composing $E$ with the restriction to the interval $I_m$.
Note that $\BV^\#(I_m)^n$ is a Banach space with the norm $\tnorm{u}_m$ defined as the sum of the variation of $u$ over $I_m$ and the norm of $u(\tau_0)$.

\begin{thm}
The operator $E_m:T_\mx\to\BV^\#(I_m)^n$ is bounded.
\end{thm}

\begin{proof}
Due to the closed graph theorem we merely have to show that $E_m$ is a closed operator.
Thus assume that the sequence $([u_j],[f_j])$ converges to $([u],[f])$ in $T_\mx$ and that $E_m([u_j],[f_j])$ converges to $v$ in $\BV^\#(I_m)^n$ and hence pointwise.
To simplify notation we assume that $E_m([u_j],[f_j]))$ and $E_m([u],[f])$ are the restrictions of $u_j$ and $u$, respectively, to the interval $I_m$.
We need to show that $u=v$ on $I_m$.

First note that $u_j^\pm(\tau_{\pm k})\in N_{\pm k}^{\perp}$ and $\abs{u_j^\pm(\tau_{\pm k})-v^\pm(\tau_{\pm k})}\rightarrow0$ imply that $v$ satisfies conditions $(\pm k)$ for each $k\in\{0,..., m-1\}$.
For $\ell\in\{-m,m-1\}$ and $x\in (\tau_\ell,\tau_{\ell+1})$ we have
$$u_j^-(x)=U_\ell^-(x)\Big(u_j^+(\tau_\ell)+J^{-1}\int_{(\tau_\ell,x)}U_\ell^*wf_j\Big)$$
when $U_\ell$ denotes the fundamental matrix of $Ju'+qu=0$ on the interval $(\tau_\ell,\tau_{\ell+1})$ satisfying $U_\ell^+(\tau_\ell)=\id$.
Taking the limit as $j\to\infty$ gives
$$v^-(x)=U_\ell^-(x)\Big(v^+(\tau_\ell)+J^{-1}\int_{(\tau_\ell,x)}U_\ell^*wf\Big)$$
since the integral may be considered as a vector of scalar products which are, of course, continuous.
The variation of constants formula shows that $v$ is a balanced solution for $Jv'+qv=wf$ on $(\tau_\ell,\tau_{\ell+1})$.
We also have
\begin{equation}\label{201012.1}
J(u_j^+(\tau_\ell)-u_j^-(\tau_\ell))+\Delta_q(\tau_\ell)u_j(\tau_\ell)=\Delta_w(\tau_\ell)f_j(\tau_\ell).
\end{equation}
The fact that $[f_j]$ converges to $[f]$ in $L^2(w)$ implies, on account of the positivity of $w$, that $\Delta_w(\tau_\ell)f_j(\tau_\ell)$ converges to $\Delta_w(\tau_\ell)f(\tau_\ell)$.
Therefore taking a limit in \eqref{201012.1} shows, in conjunction with the previous observations, that $Jv'+qv=wf$ on the interval $I_m$.
Since $u$ satisfies the same equation we have that $u-v$ satisfies $J(u-v)'+q(u-v)=0$ on $I_m$.

Next we show $w(u-v)=0$ on $I_m$.
Fatou's lemma implies
$$0\leq \int_{I_m}(u-v)^*w(u-v)\leq \liminf_{j\to\infty}\int_{I_m}(u-u_j)^*w(u-u_j)=0.$$
It follows that $w(u-v)=0$ on $I_m$.

Finally, a variant of Lemma \ref {L:5.1} shows now that $u=v$.
\end{proof}

\section{Green's function}\label{S:G}
Now suppose that we have a self-adjoint restriction $T$ of $T_\mx$.
The resolvent set of $T$ is the set of those $\la$ for which $T-\la:\dom(T)\to L^2(w)$ is bijective, i.e.,
$$\varrho(T)=\{\la\in\bb C: \ker(T-\la)=\{0\}, \ran(T-\la)=L^2(w)\}$$
which is an open set.
We denote its complement, the spectrum of $T$, by $\sigma(T)$.
Since $T$ is self-adjoint, $\sigma(T)$ is a subset of $\bb R$.
If $\la\in\varrho(T)$, then the resolvent $R_\la=(T-\la)^{-1}$ is a bounded linear operator from $L^2(w)$ to $\dom (T)$.
We now define $\mc R_\la: L^2(w)\to\BVl^\#((a,b))^n$ by
$$\mc R_\la[f]=E((R_\la[f],\la R_\la[f]+[f])).$$
Thus $\mc R_\la[f]$ is the unique solution of $Ju'+qu=w(\la u+f)$ in $\mc L^2(w)$ satisfying condition $(k)$ for every $k\in\bb Z$.

We will now show that $\mc R_\la$ is an integral operator.
Its kernel $G$ is called a Green's function for $T$.
\begin{thm}
If $T$ is a self-adjoint restriction of $T_\mx$, then there exists, for given $x\in(a,b)$ and $\la\in\varrho(T)$, a matrix $G(x,\cdot,\la)$ such that the columns of $G(x,\cdot,\la)^*$ are in $L^2(w)$ and
\begin{equation}\label{201015.1}
(\mc R_\la[f])(x)=\int G(x,\cdot,\la)wf.
\end{equation}
\end{thm}

\begin{proof}
Fix $x\in I_m$ and $\la\in\varrho(T)$.
Consider the restriction of $\mathcal{R}_\lambda[f]$ to the interval $I_m$.
Since $E_m$ and $R_\la$ are bounded operators the map $[f]\mapsto (\mathcal{R}_\lambda[f])(x)$ is a bounded linear map from $L^2(w)$ to $\bb C^n$.
Hence there are elements $[g_1], ..., [g_n]\in L^2(w)$ such that the $k$-th component of $(\mathcal{R}_\lambda[f])(x)$ equals $\langle [g_k],[f]\rangle$.
Let these be the columns of the matrix-valued function $G(x,\cdot,\lambda)^*$.
Then we obtain \eqref{201015.1}.
\end{proof}

One wishes to complement this fairly abstract existence result by a more concrete one where Green's function is given in terms of solutions of the differential equation as is done in the classical case, see, for instance, Zettl \cite{MR2170950}.
This was also achieved in \cite{MR4047968} under the assumption that $\Xi_0$ is empty and minor generalizations of this are certainly possible.
Such an explicit construction of Green's function, where possible, is the cornerstone of many other results in spectral theory, in particular the development of a spectral transformation and more detailed information about the resolvent, e.g., the compactness of the resolvent in the regular case.
Due to the difficulties posed by the absence of an existence and uniqueness theorem for initial value problems we have, so far, not been able to obtain such a construction in general.
However, we hope to return to this issue in the future.

\section{Example}
In this section we treat an example where the matrices $B_\pm(x,\la)$ fail to be invertible for infinitely many $x$ and all $\la$, in other words where $\Xi_0$ is infinite and $\Lambda_x=\bb C$ for all $x\in\Xi_0$ (recall that in \cite{MR4047968} the hypothesis $\Xi_0=\emptyset$ was made causing each $\Lambda_x$ to be finite).
The example is
$Ju'+qu=wf$ on $(a,b)=\bb R$ where
$$J=\begin{pmatrix}0&-1\\ 1&0\end{pmatrix},\;
 q=\begin{pmatrix}0&2\\ 2&0\end{pmatrix}\sum_{k\in\bb Z} (\delta_{2k}-\delta_{2k+1}),\;\text{and,}\;
 w=\begin{pmatrix}2&0\\ 0&0\end{pmatrix}\sum_{k\in\bb Z}\delta_k$$
with $\delta_k$ denoting the Dirac point measure concentrated on $\{k\}$.
Since we are seeking balanced solutions we need the matrices
$$B_-(2k-1,\la)=\begin{pmatrix}\la&0\\ 2&0\end{pmatrix}\quad\text{and}\quad B_+(2k-1,\la)=\begin{pmatrix}-\la&-2\\ 0&0\end{pmatrix}$$
as well as
$$B_-(2k,\la)=\begin{pmatrix}\la&-2\\ 0&0\end{pmatrix}\quad\text{and}\quad B_+(2k,\la)=\begin{pmatrix}-\la&0\\ 2&0\end{pmatrix}.$$
If $x$ is not an integer we have $B_\pm(x,\la)=J$.
Note that $f\in\mc L^2(w)$ if and only if $k\mapsto f_1(k)$ is in $\ell^2(\bb Z)$ and any element in $L^2(w)$ is uniquely determined by these values (here $f_1$ denotes the first component of $f$).

In any interval $(k,k+1)$ solutions of $Ju'+qu=w(\la u+f)$ are constant, say $(\alpha_k,\beta_k)^\top$.
At $x=2k-1$ the equation
$$B_+(2k-1,\la)u^+(2k-1)-B_-(2k-1,\la)u^-(2k-1)=(2f_1(2k-1),0)^\top$$
implies $\alpha_{2k-2}=0$ and
\begin{equation}\label{ev}
-\la \alpha_{2k-1}-2\beta_{2k-1}=2f_1(2k-1).
\end{equation}
Similarly, at $x=2k$ we get $\alpha_{2k}=0$ and
\begin{equation}\label{od}
-\la\alpha_{2k-1}+2\beta_{2k-1}=2f_1(2k).
\end{equation}
We can now describe the space $\mc T_\mx$.
A pair $(u,f)$ is in $\mc T_\mx$ if and only if the sequences $k\mapsto f_1(k)$ and $k\mapsto u_1(k)$ are in $\ell^2(\bb Z)$, $f_1(2k)=-f_1(2k-1)$, $u_1(2k)=u_1(2k-1)$, and
$$u=\sum_{k\in\bb Z}\Big(\bsm{2u_1(2k)\\ f_1(2k)}\chi^\#_{(2k-1,2k)} + \bsm{0\\ \beta_{2k}}\chi^\#_{(2k,2k+1)}\Big)$$
with arbitrary numbers $\beta_{2k}$.
Note that $\|u\|^2=4 \sum_{k\in\bb Z} |u_1(2k)|^2$.

Choosing here $f=0$ shows that $0$ is an eigenvalue of $T_\mx$ with infinite multiplicity.
Choosing $f=0$ and requiring $\|u\|=0$ determines the space $\mc L_0$.
Indeed,
$$\mc L_0=\Big\{\sum_{k\in\bb Z} \bsm{0\\ \beta_{2k}} \chi^\#_{(2k,2k+1)}: \beta_{2k}\in\bb C\Big\}$$
which is infinite-dimensional.
We now define the sequence $\tau$ setting $\tau_0=1/2$ and, for $k\in\bb N$, $\tau_k=k$ and $\tau_{-k}=1-k$.
A solution $u$ of $Ju'+qu=w(\la u+f)$ always satisfies condition $(2k+1)$ and it satisfies condition $(2k)$ exactly when $\beta_{2k}=0$.

For $f=0$ equations \eqref{ev} and \eqref{od} show that no non-zero $\la$ can be an eigenvalue of $T_\mx$.
In particular, the deficiency indices $n_\pm$ are $0$, i.e., $T_\mx$ is self-adjoint.
Now choose $\la\neq0$ and $f$ arbitrary in $L^2(w)$.
Then
\begin{equation}\label{201015.2}
(\mc R_\la f)(x)=-\frac{1}{2\la} \sum_{k\in\bb Z}\begin{pmatrix}2f_1(2k-1)+2f_1(2k)\\ \la f_1(2k-1)-\la f_1(2k)\end{pmatrix}\chi^\#_{(2k-1,2k)}(x)
\end{equation}
is the unique solution of $Ju'+qu=w(\la u+f)$ satisfying condition $(k)$ for any $k\in\bb Z$.
Since
\begin{equation}\label{nrla}
\|\mc R_\la f\|^2=\sum_{k\in\bb Z} 2|(\mc R_\la f)_1(k)|^2=\frac1{|\la|^2} \sum_{k\in\bb Z} |f_1(2k-1)+f_1(2k)|^2
\end{equation}
is finite we have that $\bb C\setminus\{0\}$ is the resolvent set of $T_\mx$.

We now define $\mc H=\{u\in L^2(w): u_1(2k-1)=u_1(2k)\}$ and $\mc H_\infty=\{f\in L^2(w): f_1(2k-1)=-f_1(2k)\}$.
These spaces are orthogonal to each other and their direct sum is $L^2(w)$.
Equation \eqref{nrla} shows that $\ker R_\la=\mc H_\infty$.
Moreover, we have
$$T_\mx=(\mc H\times\{0\})\oplus (\{0\}\times\mc H_\infty).$$
This is an instance of a general feature for a self-adjoint linear relation $T$: if $\mc H$ is the closure of the domain of $T$, $\mc H_\infty$ the orthogonal complement of $\mc H$, and $T_0=T\cap(\mc H\times\mc H)$, then $T=T_0\oplus (\{0\}\times\mc H_\infty)$.
The former summand is then a linear operator densely defined in $\mc H$ called the operator part of $T$.
The latter summand is called the multi-valued part of $T$.

We end this example by identifying Green's function for our example.
It may be guessed by looking at equation \eqref{201015.2}.
In any case one can check directly that $(\mc R_\la f)(x)=\int G(x,\cdot,\la)wf$.
Note that the second column of $G$ is irrelevant since the second row of $w$ is $0$.
When $x$ is not integer $G(x,y,\la)$ is given by
$$\sum_{k\in\bb Z}\left[-\frac1{\la}\begin{pmatrix}1&0\\ 0&0\end{pmatrix}+\frac12\begin{pmatrix}0&1\\ -1&0\end{pmatrix}\sgn(x-y)\right] \chi^\#_{(2k-1,2k)}(x)\chi^\#_{(2k-1,2k)}(y).$$
If $x$ is an integer we have instead
$$G(2k-1,y,\la)=\frac12 \lim_{x\downarrow 2k-1} G(x,y,\la) \quad\text{and}\quad G(2k,y,\la)=\frac12 \lim_{x\uparrow 2k} G(x,y,\la).$$


\end{document}